\documentclass[12pt]{article}
\usepackage{fancyhdr, array,calc,graphicx,url,tabularx}
\usepackage{geometry}
\usepackage{latexsym}
\usepackage{amssymb}
\usepackage{amsmath}
\usepackage{makeidx}
\usepackage{enumerate}
\usepackage{verbatim}
\usepackage{tikz}
\usepackage[colorlinks=true,linkcolor=blue]{hyperref}

\usepackage{changepage}

\makeindex

{
   \end{minipage}
   \vspace*{\stretch{3}}
   \clearpage
}

\def\undertilde#1{{\baselineskip=0pt\vtop
  {\hbox{$#1$}\hbox{$\scriptscriptstyle\sim$}}}{}}

\newcommand{\utilde}{\undertilde}

\renewcommand{\gg}{\gamma}

\newcommand{\bR}{{\mathbb{R}}}

\newcommand{\rest}{\restriction}

\renewcommand{\models}{\vDash}
\newcommand{\powerset}{{\wp}}
\newcommand{\cp}{{\rm crit }}

\newcommand{\lh}{{\rm lh}}

\newtheorem{theorem}{Theorem}[section]

\newtheorem{lemma}[theorem]{Lemma}
\newtheorem{corollary}[theorem]{Corollary}

\numberwithin{figure}{section}

\newenvironment{proof}{{\it{
Proof.}}}{\nopagebreak\mbox{}{\hfill$\square$}
\par\bigskip}
\newenvironment{Claim}{\noindent{\it
Claim} {}}{}
\newenvironment{Case}{\noindent{\it
Case}}{}

\newcommand{\rthm}[1]{Theorem~\ref{#1}}
\newcommand{\rlem}[1]{Lemma~\ref{#1}}

\newcommand{\rcor}[1]{Corollary~\ref{#1}}

\def\inseg{\trianglelefteq}

\def\k{\kappa}
\def\a{\alpha}
\def\b{\beta}
\def\d{\delta}

\def\l{\lambda}

\def\P{{\mathcal{P} }}
\def\W{{\mathcal{W} }}
\def\Q{{\mathcal{ Q}}}
\def\mH{{\mathcal{ H}}}
\def\K{{\mathcal{ K}}}

\def\R{{\mathcal R}}
\def\X{{\mathbb X}}
\def\H{{\rm{HOD}}}
\def\M{{\mathcal{M}}}
\def\N{{\mathcal{N}}}

\def\T {{\mathcal{T}}}
\def\U{{\mathcal{U}}}
\def\S{{\mathcal{S}}}

\def\X{{\mathcal{X}}}
\def\Y{{\mathcal{Y}}}

\def\and{\mathrel{\kern1pt\&\kern1pt}}

\def\inseg{\triangleleft}
\def\insegeq{\trianglelefteq}

\def\<#1>{\langle\,#1\,\rangle}

 \input xy
 \xyoption{all}
       
\title{A characterization of extenders of $\H$\thanks{2000 Mathematics Subject Classifications:
03E15, 03E45, 03E60.}
\thanks{Keywords: Mouse, inner model theory, descriptive set theory, hod mouse.}
\thanks{The author's research was partially supported by the NSF Career Award DMS-1352034.}}

\author{Grigor Sargsyan}   
           
\date{\today}

\begin{document}

\maketitle
\begin{abstract}
Assume $AD+V=L(\mathbb{R})$. Let $\kappa=\utilde{\delta}^2_1$, the supremum of all $\utilde{\Delta}^2_1$ prewellorderings. We prove that extenders on the sequence of $\H$ that have critical point $\kappa$ are generated by countably complete measures. This provides a partial reversal of Woodin's  result that the $<\Theta$-strongness of  $\kappa$ in $\H$ is witnessed by $\kappa$-complete ultrafilters on $\k$. The aforementioned characterization of extenders works in a more general setting for all cutpoint measurable cardinals of $\H$ in all models of determinacy where the fine structural analysis of $\H$ has been carried out. For example, it holds in the minimal model of the Largest Suslin Axiom.  It also gives a simple proof of a theorem of Steel that the successor members of the Solovay sequence are cutpoints in $\H$ (in models where $\H$ analysis is carried out).
\end{abstract}

\section{The context}

During the past 20-30 years, largely because of the pioneering work of Steel and Woodin, the models of determinacy have been extensively studied  using methods from \textit{inner model theory}, which is the study of $L$-like models of fragments of set theory. 
$\H$ of a model of determinacy has been isolated as a key model to analyze. Recall that $\H$ is the inner model consisting of hereditarily ordinal definable sets. It has been shown that $\H$'s of many determinacy models are $L$-like models satisfying $GCH$, $\square$ and etc. 

The reader can learn more about the aforementioned work by consulting \cite{BSL}, \cite{Steel1995}, \cite[Chapter 8]{OIMT},\cite{HODCoreModel} or \rthm{hod of lr}.

The canonical models studied by inner model theorists have the form $L[\vec{E}]$ where $\vec{E}$ is a \textit{sequence of extenders}. Most standard textbooks of set theory treat extenders (for example, see \cite{Jech}). Perhaps the easiest way of introducing extenders is via the associated elementary embeddings.

Suppose $M$ and $N$ are two transitive models of set theory and $j:M\rightarrow N$ is a non-trivial elementary embedding. Let $\kappa=\cp(j)$ and let $\l\in [\k, j(\k))$ be any ordinal. Set 
\begin{center}
$E_j=\{(a, A)\in [\l]^{<\omega}\times \powerset(\k)^M: a\in j(A)\}$.
\end{center}
$E_j$ is called the $(\k, \l)$-extender derived from $j$. $E_j$ is really an $M$-extender as it measures the sets in $M$. As with more familiar ultrafilters, one can define extenders abstractly without using the parent embedding $j$, and then show that each extender, via an ultrapower construction, gives rise to an embedding. Given a $(\k, \l)$-extender $E$ over $M$, we let $\pi_E: M\rightarrow Ult(M, E)$ be the ultrapower embedding.  A computation that involves chasing the definitions shows that $E$ is the extender derived from $\pi_E$. Similar computations also show that $\k=\cp(\pi_E)$ and $\pi_E(\k)\geq \l$. It is customary to write $\cp(E)$ for $\k$ and $lh(E)=\l$\footnote{``$lh(E)$ is the length of $E$".}.

 The motivation behind extenders is the fact that extenders capture more of the universe in the ultrapower than one can achieve via the usual ultrapower construction. In particular, under large cardinal assumptions, one can have $(\k, \l)$-extender $E$ such that $V_\l\subseteq Ult(V, E)$. Because of this all large cardinal notions below superstrong cardinals can be captured by extenders.
 
 The extenders as we defined them above are sometimes called \textit{short} extenders, where shortness refers to the fact that all the measures of the extender concentrate on its critical point. Large cardinal notions such as supercompactness, hugeness and etc cannot be captured by such short extenders as embeddings witnessing supercompactness gives rise to measures that do not concentrate on the critical point of the embedding. However, one can capture these large cardinal notions by using the so-called \textit{long} extenders. We do not need them in this paper, and so we will not dwell on them. It is, however, important to keep in mind that in this paper extenders will always be short implying that the large cardinal notions that we are concerned with are all below superstrong cardinals. This restriction is necessary as the theory of extender models with long extenders has not yet been fully developed. 

Most of inner model theory is motivated by the inner model program, which is the program of building $L$-like models that have large cardinals. As all large cardinals postulate the existence of extenders, it is natural to look for such canonical models among the models of the form $L[\vec{E}]$ where $\vec{E}$ is an extender sequence. One wants that such models be somehow nicely definable and have $L$-like properties. Because of this, the extender sequence $\vec{E}$ has to satisfy stringent requirements, as otherwise any set can be coded into $L[\vec{E}]$ for some $\vec{E}$. Readers interested in learning the exact definition of a premouse should consult \cite[Definition 2.4]{OIMT}.

Models of the form $L_\a[\vec{E}]$ that have the desired definability and $L$-likeness conditions are called \textit{pre-mice}. Mice are \textit{iterable} premice.

 It takes some time to define iterability, but a reader familiar with iterated ultrapower construction can develop a good grasp of what this might mean. Recall that a given a $\k$-complete normal ultrafilter $\mu$ on $\k$, we can iteratively build ultrapowers of $V$. More precisely, given an ordinal $\nu$, we build a sequence $(M_\a,  j_{\a, \b}: \a<\b< \nu)$ such that
 \begin{enumerate}
 \item $M_0=V$,
 \item $j_{\a, \b}:M_\a\rightarrow M_\b$,
 \item for $\a<\b<\gg$, $j_{\a, \gg}=j_{\b, \gg}\circ j_{\a, \b}$,
 \item $M_{\a+1}=Ult(M_\a, j_{0, \a}(\mu))$ and $j_{\a, \a+1}$ is the ultrapower embedding,
 \item for a limit ordinal $\l<\nu$, $M_\l$ is the direct limit of $(M_\a, j_{\a, \b}: \a<\b<\l)$ and $j_{\a, \l}:M_\a\rightarrow M_\l$ is the direct limit embedding.
 \end{enumerate}
It is a celebrated theorem of Gaifman (see \cite[Theorem 19.7]{Jech}) that such an iteration only produces well-founded models. 

When the universe has larger cardinals, it is possible to produce iterations that do not have the above linear form but rather resemble a tree. For example, it is possible to have an iteration in which there is a model $M_\a$ and an extender $E_\a\in M_\a$ such that for some $\b<\a$, it makes sense to build $Ult(M_\b, E_\a)$\footnote{This, for instance, can happen when $\powerset(\cp(E_\a))^{M_\b}=\powerset(\cp(E_\a))^{M_\a}$.}. The resulting iteration, then, will have a tree structure. For such iterations, it is not clear how to construct a direct limit at a limit stage as there can be many branches. We then say that $\M=L_\a[\vec{E}]$ is iterable when there is a strategy of picking a branch at limit stages such that the direct limit along this branch is well-founded. The key here is that the branch is picked according to a strategy. To learn more about iteration trees one can consult \cite[Chapter 3]{OIMT}.

The references of \cite{OIMT} contain papers that establish connections with inner model theory and topics in descriptive set theory. In this paper, we will mainly deal with the following central theorem. Recall that $\Theta$ is the least ordinal that is not a surjective image of $\bR$.

\begin{theorem}[Steel, Woodin, \cite{HODCoreModel}]\label{hod of lr} Assume $AD+V=L(\bR)$. Then $V_\Theta^\H$ is a universe of a premouse, i.e., there is a premouse of the form $L_\Theta[\vec{E}]$ such that 
\begin{center}
$(V_\Theta^\H, \vec{E}, \in)=L_\Theta[\vec{E}]$.
\end{center}
\end{theorem}

One can show that the premouse representation of  $V_\Theta^\H$ is unique (for example, via the methods of \cite{MeasuresInMice}).

Working in $L(\bR)$, let $\k$ be the supremum of all $\utilde{\Delta}^2_1$-prewellorderings of $\bR$. Another theorem that is very relevant to us is the following.

\begin{theorem}[Woodin]\label{woodin} Assume $AD+V=L(\bR)$ (or $AD^++V=L(\powerset(\bR))$). Then the following statements hold true.
\begin{enumerate}
\item $\H\models ``\Theta$ is a Woodin cardinal".
\item $\k$ is the least cardinal that is $<\Theta$-strong in $V_\Theta^\H$.
\item  For each $\l\in (\k, \Theta)$, there is a $\k$-complete normal ultrafilter $\mu$ on $\k$ such that $\pi_\mu(\H)|\l=\H|\l$. 
\end{enumerate}
\end{theorem}

Clause 1 of \rthm{woodin} is \cite[Theorem 4.14]{KoelWoodin}, Clause 2 for $L(\bR)$ is shown in \cite[Chapter 8]{OIMT} and Clause 3 is \cite[Theorem 4.13]{KoelWoodin}\footnote{It is the proof that gives this result.}.  

Part 3 of the above theorem says that the extenders witnessing $<\Theta$-strongness of $\k$ in $\H$ come from ultrafilters on $\k$. Let $L_\Theta[\vec{E}]$ be the mouse representation of $V_\Theta^\H$. Part 3 of the theorem does not, however, show that all extenders in $\vec{E}$ whose critical point is $\k$ arise from ultrafilters on $\k$. \rthm{main theorem} shows that the reverse is also true for many extenders.

Given a premouse (or any model with an extender sequence) $\M$ and an $\M$-cardinal $\nu$ we let $o^\M(\nu)=\sup(\{lh(E): E\in \vec{E}^\M \wedge \cp(E)=\nu\})$. $o^\M(\nu)$ is the Mitchell order of $\nu$.

\begin{theorem}\label{main theorem} Assume $AD+V=L(\bR)$. Let $\kappa$ be the supremum of $\utilde{\Delta}^2_1$-prewellorderings of $\bR$ and let $\mathcal{H}$ be the premouse representation of $V_\Theta^\H$. Fix $E\in \vec{E}^\mathcal{H}$  such that $\cp(E)=\k$ and $\mH\models ``o^{Ult(\M, E)}(\k)$ is an inaccessible cardinal". Then there is a countably complete measure $U_E$ on $\kappa$ and an embedding $i: Ult(\mH, E)\rightarrow \pi_{U_E}(\mH)$ such that 
\begin{center}
$\pi_{U_E}\rest \mH=i\circ \pi_E$, 
\end{center}
where $\pi_E:\mH\rightarrow Ult(\mH, E)$ is the ultrapower embedding. In particular, $\pi_E(\k)\leq \pi_{U_E}(\k)$.
\end{theorem}

The inequality $\pi_E(\k)\leq \pi_{U_E}(\k)$ is rather useful, and we expect that it will have many applications. First, we note that \rthm{main theorem} is not specific to $L(\bR)$. It holds in any model $M$ of $AD^+$ that permits $\H$ analysis. Here we say that a transitive model $M\models AD^+$ permits $\H$ analysis if $V_\Theta^{\H^M}$ has a fine structural representation.   Experts of descriptive inner model theory know that rank initial segments of $\H^M$ do not in general have a representation as a premouse. Nevertheless, many have been shown to be \textit{hod premice}, which are structures of the form $L_\a[\vec{E}, \Sigma]$ where $\vec{E}$ is an extender sequence as before and $\Sigma$ is an iteration strategy for $L_\a[\vec{E}, \Sigma]$ itself. A reader who would like to learn more about $\H$ analysis can consult \cite{BSL}, \cite{ATHM}, \cite{HODCoreModel} or \cite{NormComp}.

Also, \rthm{main theorem} is not specific to $\kappa$ as defined in \rthm{main theorem}. It holds for any \textit{cutpoint} cardinal of $\H$. Here, we say that a cardinal $\nu$ is a cutpoint of a premouse $\M$ if there is no extender $E\in \vec{E}^\M$ such that $\nu\in (\cp(E), \lh(E))$.

Perhaps the clearest way to say that a model of $AD^+$ permits hod analysis is to say that it satisfies \textit{Generation}. For technical reasons, it is easier to work with $AD_{\mathbb{R}}$. Assume $AD_{\mathbb{R}}$. We say Generation holds if for any set of reals $A\in \powerset(\bR)\cap M$, there is a hod pair $(\P, \Sigma)$\footnote{We say that $(\P, \Sigma)$ is a hod pair if $\P$ is a hod mouse and $\Sigma$ is its strategy. It is usually required that the strategy has nice properties. See  \cite{ATHM} or \cite{NormComp}.} such that $A\leq_w Code(\Sigma)$\footnote{$Code(\Sigma)$ is the set of reals coding $\Sigma$, and $<_w$ is the Wadge order.}.  Generation for a collection of sets of reals rather than for sets of reals first appeared in \cite{ATHM} where it was called Generation of Full Pointclass instead of Generation (of sets of reals) for technical reasons having to do with the lack of a general comparison theorem for hod pairs. In \cite{ATHM}, the notion of hod mouse was only developed below a theory $AD_{\mathbb{R}}+``\Theta$ is regular" and in \cite{LSA}, it was developed further up to the minimal model of the Largest Suslin Axiom. Prior to \cite{ATHM}, Steel and Woodin had introduced \textit{capturing} notions\footnote{In set theory, one usually uses the word ``capturing" to indicate that the notion in question can be characterized by canonical models of fragments of set theory. For example, one says that $\Sigma^1_2$-definability is captured by the constructible universe.} such as \textit{Mouse Capturing} (e.g., \cite{DMATM})\footnote{In the author's view, capturing notions, while similar to generation and usually imply generation, have a somewhat different flavor.}. In seminal \cite{NormComp}, Steel developed the notion of a hod mouse for short extenders and introduced \textit{Hod Pair Capturing} (HPC), which is another way of stating our Generation (see the comment after \cite[Definition 0.7]{NormComp}). 

\begin{theorem}[Steel, \cite{NormComp}] Assume $AD_{\mathbb{R}}+V=L(\powerset(\bR))$ and that HPC holds. Then $V_\Theta^\H$ has a representation as a hod premouse.
\end{theorem}

\rthm{main theorem} is then true in any model of $AD_{\mathbb{R}}+V=L(\powerset(\bR))$ that satisfies HPC.

\begin{theorem}\label{main theorem 1} Assume $AD_{\mathbb{R}}+V=L(\powerset(\bR))$ and that HPC holds, and let $\mH$ be the hod premouse representation of $V_\Theta^\H$. Suppose $E\in \vec{E}^{\mH}$ is such that 
\begin{enumerate}
\item $E$ is total\footnote{I.e. measures all subsets of $\cp(E)$ in $\mH$.} over $\mH$ and
\item letting $\k=\cp(E)$, $\mH\models ``\k$ is a cutpoint".
\end{enumerate}
Then there is a countably complete measure $U_E$ on $\kappa$ and an embedding $i: Ult(\mH, E)\rightarrow \pi_{U_E}(\mH)$ such that 
\begin{center}
$\pi_{U_E}\rest \mH=i\circ \pi_E$, 
\end{center}
where $\pi_E:\mH\rightarrow Ult(\mH, E)$ is the ultrapower embedding. In particular, $\pi_E(\k)\leq \pi_{U_E}(\k)$.
\end{theorem}

The motivation behind proving \rthm{main theorem 1} is to use the condition $\pi_E(\k)\leq \pi_{U_E}(\k)$ to establish bounds on the height of the Mitchell order of $\k$. Recalling Kunen's theorem, namely that under determinacy all countably complete ultrafilters are ordinal definable,  we immediately get that

\begin{corollary}\label{corollary} $o^\M(\k)<\theta_{\powerset(\k)}$\footnote{Recall that $\theta_X$ is the supremum of all $\a$ such that there is an $OD(X)$ surjection $f:X\rightarrow \a$.}. 
\end{corollary}

The proof of \rthm{cut point suslins} contains the proof of \rcor{corollary}. \rthm{cut point suslins} confirms Woodins conjecture that, assuming HOD analysis is successful, all members of \textit{the Solovay sequence} are cutpoints of $\H$. Recall that under $AD$, the Solovay sequence is a closed sequence of ordinals $(\theta_\a: \a\leq \Omega)$ such that
\begin{enumerate}
\item $\theta_0=\theta_{\bR}$,
\item $\theta_{\a+1}=\theta_{A}$ where $A\subseteq \bR$ is any set of Wadge rank $\theta_\a$, and
\item for a limit ordinal $\l$, $\theta_\l=\sup_{\a<\l}\theta_\a$.
\end{enumerate}

\begin{theorem}[Steel]\label{cut point suslins} Assume $AD^++V=L(\powerset(\bR))$ and suppose HPC holds. Then every member of the Solovay sequence is a cutpoint cardinal of $\H$.
\end{theorem}
\begin{proof} The proof presented below is due to the author. 
It is enough to show that each $\theta_{\a+1}$ is a cutpoint in $\H$ (for $\a\geq -1)$. Towards a contradiction assume that $\theta_{\a+1}$ is not a cutpoint in $\mH$, the hod mouse representation of $V_\Theta^{\H}$. Let $\k<\theta_{\a+1}$ be the least $<\theta_{\a+1}$-strong cardinal of $\mH$. As $\mH\models ``\theta_{\a+1}$ is a Woodin cardinal", we have that $\mH\models ``\theta_{\a+1}$ is an inaccessible cardinal". Let $E\in \vec{E}^\mH$ be such that 
\begin{enumerate}
\item $\cp(E)=\k$,
\item $\pi_E(\k)>\theta_{\a+1}$ and
\item $\theta_{\a+1}$ is a cut point of $Ult(\mH, E)$.\footnote{It should be noted that we are assuming that Schlutzenberg's results from \cite{MeasuresInMice} carry over to our context. According to this theorem the existence of any extender $E$ witnessing that $\k$ is $\theta_{\a+1}$-strong implies the existence of such an extender that is on the extender sequence of $\H$. That this is indeed the case is something that the author has not verified in a published article, but he did carry out similar calculations and is confident that the result holds at least in the case relevant to this proof.}
\end{enumerate} 
Using \rthm{main theorem 1} fix a countably complete ultrafilter $U$ such that $\pi_E(\k)\leq \pi_{U}(\k)$.  Clause 2 above implies that $\theta_{\a+1}<\pi_{U}(\k)$.

By Kunen's theorem (see \cite[Theorem 7.6]{DMT}) $U$ is ordinal definable. Because $\pi_{U}(\k)=\sup\{ [f]_{U}: f:\k\rightarrow \k\}$, it follows from the aforementioned theorem of Kunen  that there is an ordinal definable surjection $h:\powerset(\k)\rightarrow \pi_{U}(\k)$. Fix a set of reals of Wadge rank $\theta_\a$ and an $OD(A)$ surjection $g:\bR\rightarrow \k$. It follows from the Coding Lemma that there is an $OD(A)$ surjection $k:\bR\rightarrow \powerset(\k)$. Setting $\tau=h\circ k$, we get that $\tau\in OD(A)$ and $\tau: \bR\rightarrow  \pi_{U}(\k)$ is a surjection.  Because $\theta_{\a+1}<\pi_{U}(\k)$, we get an $OD(A)$ surjection $m: \bR \rightarrow \theta_{\a+1}$, contradiction.
\end{proof}

Woodin was the first to notice that \rthm{cut point suslins} must be true. Sometime in 2006-2009 he gave a sequence of informal lectures presenting his ideas to the author and John Steel.  \rthm{cut point suslins} was first proven by Steel independently of the author in 2016. The author's motivation was to prove Generation not \rthm{cut point suslins}. Woodin's ideas for proving \rthm{cut point suslins} were very different from both the proof presented above and from the calculations carried out by Steel. 

Finally, we would like to stress that the main contribution of \rthm{main theorem 1} is the inequality $\pi_E(\k)\leq \pi_{U_E}(\k)$. It is a useful tool in showing that direct limit constructions, the kind of constructions used in $\H$ analysis, yield bounded structures. We hope that it will be an essential piece in the eventual proof of Generation from $AD^++L(\powerset(\bR))$. This hope is grounded in the fact that all the current approaches to Generation involve showing that certain direct limit constructions produce bounded below $\Theta$ structures (see \cite{GenericGenerators}). The conjecture that Generation is a theorem of $AD^++L(\powerset(\bR))$ is the most central conjecture of descriptive inner model theory (see  \cite{BSL} or the introduction of \cite{NormComp}). Our main motivation for proving \rthm{main theorem 1} is that it will be helpful in settling Generation. 

\section{The proof}

Here we present the proof of \rthm{main theorem}. The proof of \rthm{main theorem 1} is very similar and has some fine structure theory. The author also noticed the proof works for any direct limit model not just $\mH$. John Steel has recently circulated notes in which he presented a proof of \rthm{main theorem 1} (see \cite{JLS}) for general direct limits, and so the interested reader can read this more general proof there. Presenting the proof of \rthm{main theorem} from scratch, while certainly not an impossible task, requires adding many more pages to the current paper. Sadly these pages will only contain  what has already been well presented in other publications. Because of this we assume the standard terminology of $\H$ analysis as presented for example in \cite{HODCoreModel}.\\

\textbf{The proof of \rthm{main theorem}}\\

We assume that the extenders are indexed according to Jensen's indexing scheme, i.e., if $\M$ is a mouse and $E\in \vec{E}^\M$ is an extender then $\a$ is the index of $E$ if $\a$ is the $\M|\a$-successor of $\pi_E^{\M|\a}(\cp(E))$.
Fix an extender $E\in \vec{E}^{\mH}$ such that $\cp(E)=\k$. Let 
\begin{center}
$\l=o^{Ult(\mH, E)}(\k)=\sup\{ lh(F): F\in \vec{E}^{\mH|lh(E)}\wedge \cp(F)=\k\}$. 
\end{center}
We are assuming that $\mH\models ``\l$ is an inaccessible cardinal".

We will define $U_E$ by a definition that gives the reduction of $U_E$ to the Martin's Measure on Turing degrees (see \cite{Kanamori}). Recall that under $AD$, given a set $C$ of Turing degrees, there is a Turing degree $e$ such that either
\begin{enumerate}
\item for every Turing degree $d$, $e\leq_T d$ implies  $d\in C$,

or

\item for every Turing degree $d$, $e\leq_T d$ implies  $d\in C^c$.
\end{enumerate}
Martin's measure is the ultrafilter generated by \textit{cones}, i.e., sets of the form $A_e=\{d: e\leq_T d\}$. We let $\mu_M$ be the Martin's measure.

We say that $(\P, A)$ captures $E$ if
\begin{enumerate}
\item $\P$ is suitable,
\item $A$ is an ordinal definable set of reals,
\item $\P$ is $A$-iterable, 
\item $\gg_{A, \infty}>lh(E)$, and
\item $E\in rng(\pi_{(\P, A), \infty})$.
\end{enumerate}
We then let $E^\P=\pi_{(\P, A), \infty}^{-1}(E)$ and $\k^\P=\pi_{(\P, A), \infty}^{-1}(\k)$. Also, let $\l=o^{Ult(\mH, E)}(\k)$ and for $\P$ as above, let $\l^\P=\pi^{-1}_{(\P, A), \infty}(\l)$. Let $B=\{ \P: (\P, A)$ captures $E\}$ and $C$ be the set of those reals that code elements of $B$ (via some natural method of coding).

We say that a Turing degree $d$ is $E$-large if $C_d=_{def}\{y: [y]_T<_T d\} \cap C\not =\emptyset$. Given an $E$-large turing degree $d$ let $\M_d$ be the result of simultaneously comparing all $\P\in B$ that have a code in $C_d$. We let $B_d$ be the set of those $\P\in B$ that have a code in $C_d$. There are only countably many such $\P\in B_d$, and so the aforementioned comparison halts producing a countable common iterate that is a member of $B$. Thus, $\M_d\in B$ for every $E$-large $d$. 

Given an $E$-large $d$, letting $\l_d=\l^{\M_d}$ and $\N_d=Ult(\M_d, E^{\M_d})|(\l_d^+)^{Ult(\M_d, E^{\M_d})}$, set 
\begin{center}
$\k_d=\pi_{\N_d, \infty}(\k^{\M_d})$. 
\end{center}
Notice that $\N_d$ has no Woodin cardinals, and because $Ult(\M_d, E)$ is $A$-iterable, $\N_d$ is $(\omega_1, \omega_1)$-iterable via a unique iteration strategy. It follows that $\pi_{\N_d, \infty}$ makes sense. It is the direct limit embedding via the unique strategy of $\N_d$. The reader may choose to consult Chapter 8 of \cite{OIMT} where it is shown that if $\Sigma$ is the unique strategy of $\N_d$ then $\M_\infty(\N_d, \Sigma)$ is a rank initial segment of $\mH$ (for example, see \cite[Theorem 8.20]{OIMT}).

We now define an ultrafilter $U_E$ on $\k$ by setting $D\in U_E$ if and only if for a cone of $d$, $\k_d\in D$. Because $\mu_d$ is a countably complete ultrafilter, we have that $U_E$ is a countably complete ultrafilter. It remains to show that there is an embedding 
$i: Ult(\mH, E)\rightarrow \pi_{U_E}(\mH)$ such that $\pi_{U_E}\rest \mH=i\circ \pi_E$. 

We have that 
\begin{center}
$Ult(\mH, E)=\{\pi_E(f)(a): f\in \mH \wedge a\in \l^{<\omega}\}$.
\end{center}
 Suppose now that $h: \l^{<\omega}\rightarrow \pi_{U_E}(\k)^{<\omega}$ is any function. Usually one builds realization embeddings like our $i$ above using such functions $h$ by setting $i_h(x)=\pi_{U_E}(f)(a)$ where $f, a$ are chosen such that $x=\pi_E(f)(a)$. We need to define $h$ so that 
\begin{enumerate}
\item $i_h$ is well-defined (i.e., it is independent of the choice of the pair $(f, a)$) and
\item $i_h$ is elementary.
\end{enumerate}
What follows is a construction of one such function $h$. 

Fix an $a\in \l^{<\omega}$. We say that $\P\in B$ \textit{captures} $a$ if $a\in rng(\pi_{(\P, A), \infty})$. We let 
\begin{enumerate}
\item $a^\P=\pi^{-1}_{(\P, A), \infty}(a)$,
\item $B_a=\{\P\in B: (\P, A)$ captures $a\}$,
\item $C_a$ be the set of those reals that code elements of $B_a$,
\item for a Turing degree $d$, $C_{a, d}=\{y: [y]_T<_T d\} \cap C_a$
\end{enumerate}
Given $\P\in B_a$ we define a function $f^a_\P:\k\rightarrow \k$ as follows. We say $\a<\k$ is \textit{typical} if for some Turing degree $d$, $\a=\k_d$. We then set
\begin{center}
$f^a_\P(\a)=\begin{cases}
0 &: \a\ \text{is not typical}\\
\pi_{\N_d, \infty}(a^{\M_d}) &: \a\ \text{is typical}
\end{cases}$
\end{center}
Recall that $\N_d=Ult(\M_d, E^{\M_d})$.

The function $a\rightarrow [f^a_\P]_{U_E}$ is our intended candidate for $h$. However, it is not even clear that $f^a_\P$ is well-defined let alone $a\rightarrow [f^a_\P]_{U_E}$. We start by proving that $f^a_\P$ is well-defined. Below we will write $a_d$ for $a^{\M_d}$.

\begin{lemma}\label{main lemma} For each $a\in \l^{<\omega}$, there is $\P\in B_a$ such that $f^a_\P$ is well-defined. 
\end{lemma}
 \begin{proof}
 The proof is via a reflection argument. Assume that the claim is false. We want to use $\Sigma^2_1$-reflection inside $L(\bR)$ to find $\b<\k$ such that our claim is false in $L_\b(\bR)$. The sentence that we want to reflect is the following.\\\\
 $\phi:$ There is a tuple $(A, E, a)$ such that
 \begin{enumerate}
 \item $A$ is an ordinal definable set of reals,  
 \item  $E\in \vec{E}^{\mH|\gg_{A, \infty}}$ has critical point $\k=_{def}\utilde{\d}^2_1$,
 \item letting $\l=o^{Ult(\mH, E)}(\k)$, $a\in \l^{<\omega}$ and $\mH\models ``\l$ is an inaccessible cardinal"\footnote{Recall from \cite{HODCoreModel} the meaning of $\gg_{A, \infty}$. Given an $A$-iterable $\Q$, $\gg^\Q_A=\sup(Hull_1^\Q(\{\tau^\Q_A\})\cap \d^\Q)$ where $\tau^\Q_A$ is the term relation capturing $A$ over $\Q$ and $\d^\Q$ is the Woodin cardinal of $\Q$. Then $\gg_{A, \infty}=\pi_{(\Q, A), \infty}(\gg^\Q_A)$.},
\item for every $\P$ that is $A$-iterable and such that $(E, a)\in rng(\pi_{(\P, A), \infty})$, there are two Turing degrees $d_0$ and $d_1$ such that $\P$ has a code both in $\{y: [y]_T<_Td_0\}$ and in $\{y: [y]_T<_T d_1\}$, $\k_{d_0}=\k_{d_1}$ and
 \begin{center}
 $\pi_{\N_{d_0}, \infty}(a_{d_0})\not =\pi_{\N_{d_1}, \infty}(a_{d_1})$.
 \end{center}
  \end{enumerate}
 Notice that $\phi$ is a sentence. Let then $\b$ be the least such that
 \begin{enumerate}
 \item $L_\b(\bR)\models \phi$,
 \item $L_\b(\bR)\models ZF-Powerset$.
 \end{enumerate}
 
Fix $(G, F, b)\in L_\b(\bR)$ witnessing $\phi$. We do not change our notation for the sets $B$, $B_b$, $C$ and $C_{b, d}$, and we will also use $\M_d$, $\N_d$ and $\l_d$ as if we are working with the non-reflected objects. All of these objects are now defined with respect to $(G, F, b)$. 

 We now want to produce a pair $(\P, \Sigma)$ such that  
\begin{enumerate}
\item $L_{\b}(\bR)\models ``\P$ is suitable and $G$-iterable",
\item $\P\in B_b$,
\item $\Sigma$ is an $(\omega_1, \omega_1)$-iteration strategy for $\P$ that is $L_\b(\bR)$-fullness preserving,
\item $\Sigma$ \textit{respects} $G$\footnote{Given a $\Sigma$-iteration $k:\P\rightarrow \Q$, $k(\tau_G^\P)=\tau^\Q_G.$},
\item $\Sigma$ has \textit{full normalization}\footnote{Whenever $\Q$ is a $\Sigma$-iterate of $\P$ via some iteration $p$, there is a normal tree $\T$ according to $\Sigma$ whose last model is $\Q$. Moreover, if $\pi^p:\P\rightarrow \Q$ is defined then $\pi^\T$ is defined and $\pi^p=\pi^\T$. See the discussion after \cite[Remark 2.2]{NormComp} and \cite[Theorem 1.1]{MPSC}}. 
\end{enumerate}
The production of $(\P, \Sigma)$ is completely standard. The basic idea is to pick a \textit{good} pointclass $\Gamma$ such that
$(\utilde{\Sigma^2_1})^{L_\b(\bR)}\subseteq \utilde{\Delta}_\Gamma$. An example of such a $\Gamma$ is $(\Sigma^2_1)^{L_\xi(\bR)}$ where $\xi>\b$ and $\xi$ ends a weak gap\footnote{A new $\utilde{\Sigma}_1$-fact is true in $L_{\xi+1}(\bR)$.}. One then works inside some $\N^*_x$-like model for $\Gamma$ and performs a fully backgrounded construction. This construction produces the desired pair $(\P, \Sigma)$. The reader can, for example, consult \cite[Lemma 5.18]{ATHM}, \cite[Proposition 2.2]{lpR}, \cite[Lemma 2.4]{lpR} and \cite[Chapter 4, 5]{NormComp}.
 
 We now show that $f^\P_b$ is well-defined. Towards a contradiction, fix $d_0$ and $d_1$ such that $\P$ has a code both in $C_{b, d_0}$ and $C_{b, d_1}$, and 
 \begin{center}
 $\pi_{\N_{d_0}, \infty}(a^{\M_{d_0}})\not =\pi_{\N_{d_1}, \infty}(a^{\M_{d_1}})$.
 \end{center}
 
 Because $\P$ has a code both in $C_{b, d_0}$ and in $C_{b, d_1}$,  we have that both $\M_{d_0}$ and $\M_{d_1}$ are $\Sigma$-iterates of $\P$. For $i\in 2$ let $\T_i$ be the normal $\P$-to-$\M_{d_i}$ tree according to $\Sigma$.
 
  Because $(\P, G)$ captures $F$ and $b$, and $\Sigma$ respects $G$, we have that\\\\
 (1) for $i\in 2$, $F^{\M_{d_i}}, b^{\M_{d_i}} \in rng(\pi^{\T_i})$.\\\\
  For $i\in 2$, let $\xi_i$ be least such that $F^{\M_{d_i}}\in \M_{\xi_i}^{\T_i}$. Next, set $\X_i=\T_i\rest \xi_i+1$. Notice that for $i\in 2$, the generators of $\X_i$ are contained in $\l_{d_i}$. In fact\\\\
 (2) for $i\in 2$, if $H$ is an extender used in $\X_i$ then $lh(H)<\l_{d_i}$.\\\\  
 For $i\in 2$, let $\R_i=\M_{\xi_i}^{\T_i}$, and let $\Lambda_i$ be the unique strategy of $\N_{d_i}$. Thus, $\R_i$ is the last model of $\X_i$.\\

\textit{Claim 1.}  $\M_\infty(\N_{d_0}, \Lambda_0)=\M_\infty(\N_{d_1}, \Lambda_1)$. \\\\
\begin{proof}  For $i\in 2$, let $\nu_i$ be such that $\M_\infty(\N_{d_i}, \Lambda_i)=\mH|\nu_i$. Let $\tau_i$ be the cardinal predecessor of $\nu_i$ in $\mH$\footnote{Because $\N_{d_i}$ is a rank initial segment of $Ult(\M_{d_i}, F^{\M_{d_i}})$ we have that $\nu_i$ is a successor cardinal in $\mH$. This observation is true regardless weather $\l$ is inaccessible in $\mH$ or not.}. It follows $\tau_i=o^{\mH}(\k_{d_i})$. Because $\k_{d_0}=\k_{d_1}$, we get that $\tau_0=\tau_1$. Hence, $\nu_0=\nu_1$
\end{proof}

Claim 1 above implies that  $\N_{d_0}$ and $\N_{d_1}$ compare to the same mouse (recall that $\N_{d_i}=Ult(\M_{d_i}, F^{\M_{d_i}})|(\l_{d_i}^+)^{Ult(\M_{d_i}, F^{\M_{d_i}})}$\footnote{Because of our assumption that $\l$ is inaccessible in $\mH$, $\N_{d_i}|\l_{d_i}=\M_{d_i}|\l_{d_i}$.}). Let $(\Y_0, \Y_1)$ be the result of the aforementioned coiteration where $\Lambda_i$ is used to iterate $\N_{d_i}$.

As for $i\in 2$, $\N_{d_i}|\l_{d_i}$ is a rank initial segment of $\M_{d_i}$ and $\l_{d_i}$ is inaccessible in $\M_{d_i}$, we can think of $\Y_i$ as a normal tree on $\R_i$.  More precisely, there is a tree $\Y_i^+$ on $\R_i$ whose tree structure is the same as the tree structure of $\Y_i$, and the extenders used in $\Y_i^+$ are exactly those used in $\Y_i$. The only difference between $\Y_i$ and $\Y_i^+$ are the models. Moreover, for each $\gg<lh(\Y_i^+)$, $\M^{\Y_i}_\gg\insegeq \M^{\Y_i}_\gg$ and if  $\M^{\Y_i}_\gg\inseg \M^{\Y_i}_\gg$ then $\M_\gg^{\Y_i}$ is a rank initial segment of $\M^{\Y_i^+}_\gg$ (the equality $\M_\gg^{\Y_i}=\M_\gg^{\Y_i^+}$ can happen exactly when there is a drop at $\gg$).

For $i\in 2$, let $\S_i$ be the last model of $\Y_i^+$. \\

\textit{Claim 2.} $\S_0=\S_1$.\\\\
\begin{proof} Let $\N$ be the common mouse that $\N_{d_0}$ and $\N_{d_1}$ coiterate to. Thus, $\N$ is the last model of $\Y_0$ and $\Y_1$. Therefore, $\N\inseg \S_i$ for $i\in 2$. Because $\Sigma$ has full normalization, for each $i$, there is a normal tree $\K_i$ on $\P$ with last model $\S_i$, namely the full normalization of $\X_i^\frown \Y_i^+$. Another way of obtaining $\K_i$ is by noting that it is the result of comparing $\S_i$ with $\P$. Then we must have that $\S_i$ side doesn't move producing a normal tree on $\P$ whose last model is $\S_i$.

Let $\iota_i<lh(\K_i)$ be least such that $\N\insegeq \M_{\iota_i}^{\K_i}$. We claim that $\iota_i+1=lh(\K_i)$. Let $\U$ be the normal tree on $\P$ with last model $\Q$ that is build by comparing $\N$ with $\P$. As $\K_i$ is also build this way, we must have that  for $i\in 2$, $\U$ is an initial segment of $\K_i$ and $\Q=\M_{\iota_i}^{\K_i}$. Thus, $\K_0|\iota_0+1=\K_1|\iota_1+1=\U$. Set $\iota=_{def}\iota_0=\iota_1$. 

Suppose now that $\iota+1<lh(\K_0)$. Let $Z_i^*$ be the extender used in $\K_i$ at stage $\iota$. Notice that\\\\
(a)  the generators of both $\K_0$ and $\K_1$ are contained in $\N$. \\\\
(a) is a consequence of our choice of $\X_i$ and $\Y_i$. We have that for $i\in 2$,
\begin{center}
$\S_i=\{\pi^{\X_i^\frown \Y_i^+}(f)(s): s\in \N\}$.
\end{center}
Also,\\\\
(b) if $H$ is an extender used in $\U$  then $lh(H)\in \N$.\\\\
Thus, $\cp(Z_i^*)\in  \N$. We must also have that either $lh(Z_0^*)>Ord\cap \N$ or $lh(Z_1^*)>Ord\cap \N$ (as otherwise $Z_i^*\in \vec{E}^{\N}$ which would imply that $\N\not \insegeq \S_i$)\footnote{We say ``either" because it might be the case that one of the trees is being padded.}.

Suppose now that $Z_0^*$ is defined (the other case is completely symmetric). Let $Z$ be the least extender used on the main branch of $\K_0$ after stage $\iota$ such that $\cp(Z)\in \N$ and $lh(Z)>Ord\cap \N$. Let $\nu=\cp(Z)$. Because $\K_0$ is normal, we have that (for an outline of the proof see the footnote of (c))
\begin{center}
$o.t.[(\{\pi^{\K_0}(f)(s): f\in \P \wedge s\in [\nu]^{<\omega}\}\cap \N\cap Ord)]=\nu$\footnote{o.t. stands for ``order type".}.
\end{center}
However, because of (2) and (b) we have that\\\\
(c) for every $\zeta<\l_{d_0}$, $o.t.[(\{\pi^{\X_0}(f)(s): f\in \P \wedge s\in [\zeta]^{<\omega}\}\cap \l_{d_0})]> \zeta$\footnote{If there is no extender $H$ used on the main branch of $\X_0$ such that $\zeta\in [\cp(H), lh(H))$ then $\zeta\in (\{\pi^{\X_0}(f)(s): f\in \P \wedge s\in [\zeta]^{<\omega}\}\cap \l_{d_0})$. Fix now $H$ such that $\zeta\in [\cp(H), lh(H))$ and $H$ is used on the main branch of $\X_0$. Let $\W$ be the model on the main branch of $\X_0$ to which $H$ is applied. Then $\pi^{\X_0}_{\W, \R_0}(\zeta)\leq \l_{d_0}$ (as all extenders used in $\X_0$ have lengths $<\l_{d_0}$). It follows that $\pi^{\X_0}_{\W, \R_0}(\zeta)\in (\{\pi^{\X_0}(f)(s): f\in \P \wedge s\in [\zeta]^{<\omega}\}\cap \l_{d_0})$.}  and\\
(d) for every $\zeta\in \N$, $o.t.[(\{\pi^{\Y^+_0}(f)(s): f\in\R_0  \wedge s\in [\zeta]^{<\omega}\}\cap \N\cap Ord)]> \zeta$.\\\\
Because $\pi^{\K_0}=\pi^{\Y^+_0}\circ \pi^{\X_0}$, we have that 
\begin{center}
$o.t.[(\{\pi^{\K_0}(f)(s): f\in \P \wedge s\in [\nu]^{<\omega}\}\cap \N\cap Ord)]>\nu$.
\end{center}
This is a contradiction.
\end{proof}

We set $\W=\S_0=\S_1$. To finish the proof of \rlem{main lemma}, we need to show that\\\\
(*) $\pi_{\N_{d_0}, \infty}(b_{d_0})=\pi_{\N_{d_1}, \infty}(b_{d_1})$.\\\\
Recall that $b_{d}=b^{\M_d}$ is the pre-image of $b$ in $\M_d$. The equality in (*) directly contradicts clause 4 of $\phi$, which we assume to be true in $L_\b(\bR)$. To show (*), it is enough to show\\\\
 (**) $\pi^{\Y_0^+}(b_{d_0})=\pi^{\Y_1^+}(b_{d_1})$.\\\\
 To see that (**) implies (*), notice that for $i\in 2$, 
 \begin{enumerate}
 \item $\pi^{\Y_i^+}\rest \N_{d_i}=\pi^{\Y_i}$ and
 \item $\pi_{\N_{d_i}, \infty}=\pi_{\N, \infty}\circ \pi^{\Y_i}$ 
 \end{enumerate}
 It then follows that assuming (**),
\begin{equation*}
\begin{split}
\pi_{\N_{d_0}, \infty}(b_{d_0}) & =\pi_{\N, \infty}\circ \pi^{\Y_0}(b_{d_0})\\
 & = \pi_{\N, \infty}\circ \pi^{\Y_1}(b_{d_1})\\
 & =\pi_{\N_{d_1}, \infty}(b_{d_1}) 
\end{split}
\end{equation*}

We now prove that (**) holds. Notice that for $i\in 2$, 
 \begin{center}
$\pi_{(\W, G), \infty}(\pi^{\Y_i^+}(b_{d_i}))=\pi_{(\R_i, G), \infty}(b_{d_i})$.
\end{center}
Thus, (**) follows once we show that $\pi_{(\R_0, G), \infty}(b_{d_0})=\pi_{(\R_1, G), \infty}(b_{d_1})$. But because $\R_i$ is an iterate of $\P$ via $\X_i$, we have that 
\begin{equation*}
\begin{split}
\pi_{(\R_0, G), \infty}(b_{d_0}) & =\pi_{(\P, G), \infty}(b^\P)\\
 & = \pi_{(\R_1, G), \infty}(b_{d_1}).
\end{split}
\end{equation*}
This finishes the proof of \rlem{main lemma}.
\end{proof}

 Notice that if $\P$ and $\Q$ are such that $f^a_\P$ and $f^a_\Q$ are well-defined then for $U_E$-almost all $\a$, $f^a_{\P}(\a)=f^a_{\Q}(\a)$. Indeed, if $d$ is any $(E, a)$-large Turing degree such that both $\P$ and $\Q$ have a code in $C_{a, d}$, $f^a_\P(\k_d)=f^a_\Q(\k_d)$.

Define $h: [\l]^{<\omega}\rightarrow [\pi_{U_E}(\k)]^{<\omega}$ by setting $h(a)=[f^a_{\P}]_{U_E}$ where $\P$ is chosen so that $f^a_\P$ is well-defined. Define $i: Ult(\mH, E)\rightarrow \pi_{U_E}(\mH)$ by 
\begin{center}
$\pi_E^{N}(g)(a)=\pi_{U_E}(g)(h(a))$.
\end{center}
Clearly, $\pi_{U_E}\rest \mH=i\circ \pi_E$. We need to show that $i$ is elementary. Thus, the next claim finishes the proof of \rthm{main theorem}. \\

\textit{Claim 3.} $i$ is elementary.\\\\
\begin{proof}
To see this, suppose $\pi_{U_E}(\mH)\models \phi[i(x)]$ where $x=\pi_E(g)(a)$ for some $g\in \mH$ and $a\in [\l]^{<\omega}$. Let $\P\in B_a$ be such that $f^a_\P$ is well-defined. Then
\begin{center}
 $\pi_{U_E}(\mH)\models \phi[\pi_{U_E}(g)([f^a_{\P}]_{U_E})]$.
\end{center} 
For an $(E, a)$-large degree $d$, let $a_{d, \infty}=\pi_{\N_d, \infty}(a_d)$. Chasing definitions we get that
\begin{center} $\{ d:\mH\models \phi[g(a_{d, \infty})]\}\in \mu_M$.\end{center}
Fix $d$ in the above set such that $g\in rng(\pi_{(\M_d, A), \infty})$ and let $k$ be the $\pi^{-1}_{(\M_d, A), \infty}(g)$. Notice that letting $\W=Ult(\M_d, E^{\M_d})$, $a_{d, \infty}=\pi_{(\W, A), \infty}(a_d)$. It follows that 
\begin{center}
$Ult(\M_d, E^{\M_d})=\W\models  \phi[\pi_{E^{\M_d}}(k)(a_d)]$.
\end{center}
Hence, 
\begin{center} $\M_d \models \{ t:\phi[k(t)]\}\in (E^{\M_d})_{a_d}$.\end{center}
Applying $\pi_{(\M_d, A), \infty}$, we get that 
\begin{center} $\mH\models \{ t:\phi[g(t)]\}\in E_{a}$.\end{center}
Therefore, $Ult(\mH, E)\models \phi[\pi_E(g)(a)]$.
\end{proof}

\bibliographystyle{plain}
\bibliography{CharExteHod.bib}

\begin{thebibliography}{10}

\bibitem{Blah}
Andr\'{e}s~Eduardo Caicedo, Paul Larson, Grigor Sargsyan, Ralf Schindler, John
  Steel, and Martin Zeman.
\newblock Square principles in {${\Bbb P}_{\max}$} extensions.
\newblock {\em Israel J. Math.}, 217(1):231--261, 2017.

\bibitem{Jech}
Thomas Jech.
\newblock {\em Set theory}.
\newblock Springer Monographs in Mathematics. Springer-Verlag, Berlin, 2003.
\newblock The third millennium edition, revised and expanded.

\bibitem{Kanamori}
Akihiro Kanamori.
\newblock {\em The higher infinite}.
\newblock Perspectives in Mathematical Logic. Springer-Verlag, Berlin, 1994.
\newblock Large cardinals in set theory from their beginnings.

\bibitem{KoelWoodin}
Peter Koellner and W.~Hugh Woodin.
\newblock Large cardinals from determinacy.
\newblock In {\em Handbook of set theory. {V}ols. 1, 2, 3}, pages 1951--2119.
  Springer, Dordrecht, 2010.

\bibitem{GenericGenerators}
Grigor Sargsyan.
\newblock Generic generators.
\newblock {\em To appear.}

\bibitem{BSL}
Grigor Sargsyan.
\newblock A short tale of hybrid mice, to appear at the {B}ulletin of
  {S}ymbolic {L}ogic.

\bibitem{ATHM}
Grigor Sargsyan.
\newblock {\em Hod mice and the mouse set conjecture}, volume 236 of {\em
  Memoirs of the {A}merican {M}athematical {S}ociety}.
\newblock American Mathematical Society, 2014.

\bibitem{lpR}
Grigor Sargsyan and John Steel.
\newblock The mouse set conjecture for sets of reals.
\newblock {\em J. Symb. Log.}, 80(2):671--683, 2015.

\bibitem{LSA}
Grigor Sargsyan and Nam Trang.
\newblock The {L}argest {S}uslin {A}xiom.
\newblock 2013.
\newblock available at math.rutgers.edu/$\sim$gs481/books.

\bibitem{MeasuresInMice}
Farmer~Salamander Schlutzenberg.
\newblock {\em Measures in mice}.
\newblock ProQuest LLC, Ann Arbor, MI, 2007.
\newblock Thesis (Ph.D.)--University of California, Berkeley.

\bibitem{MPSC}
John~R. Steel.
\newblock Mouse pairs and {S}uslin cardinals.
\newblock {\em To appear, available at https://math.berkeley.edu/$\sim$steel/}.

\bibitem{NormComp}
John~R. Steel.
\newblock Normalizing iteration trees and comparing iteration strategies.
\newblock {\em To appear, available at https://math.berkeley.edu/$\sim$steel/}.

\bibitem{JLS}
John~R. Steel.
\newblock Notes on work of {J}ackson and {S}argsyan.
\newblock {\em available at https://math.berkeley.edu/$\sim$steel/}.

\bibitem{Steel1995}
John~R. Steel.
\newblock {${\rm HOD}\sp {{\rm{L}}(\Bbb R)}$} is a core model below {$\Theta$}.
\newblock {\em Bull. Symbolic Logic}, 1(1):75--84, 1995.

\bibitem{DMATM}
John~R. Steel.
\newblock Derived models associated to mice.
\newblock In {\em Computational prospects of infinity. {P}art {I}.
  {T}utorials}, volume~14 of {\em Lect. Notes Ser. Inst. Math. Sci. Natl. Univ.
  Singap.}, pages 105--193. World Sci. Publ., Hackensack, NJ, 2008.

\bibitem{DMT}
John~R. Steel.
\newblock The derived model theorem.
\newblock In {\em Logic {C}olloquium 2006}, Lect. Notes Log., pages 280--327.
  Assoc. Symbol. Logic, Chicago, IL, 2009.

\bibitem{OIMT}
John~R. Steel.
\newblock An outline of inner model theory.
\newblock In {\em Handbook of set theory. {V}ols. 1, 2, 3}, pages 1595--1684.
  Springer, Dordrecht, 2010.

\bibitem{HODCoreModel}
John~R. Steel and W.~Hugh Woodin.
\newblock {H}{O}{D} as a core model.
\newblock In {\em Ordinal definability and recursion theory: {T}he {C}abal
  {S}eminar. {V}ol. {III}}, volume~43 of {\em Lect. Notes Log.}, pages
  257--345. Assoc. Symbol. Logic, Ithaca, NY, 2016.

\end{thebibliography}
\end{document}